\documentclass[11pt]{article}

\usepackage{graphicx}
\usepackage{amsmath,amsfonts,amssymb}
\usepackage{citesort}
\usepackage{url}
\usepackage{amsmath}
\usepackage{graphicx}
\usepackage{epsfig}
\usepackage{epstopdf}
\usepackage{color}
\def\tto{\;{\lower 1pt \hbox{$\rightarrow$}}\kern -10pt
\hbox{\raise 2pt \hbox{$\rightarrow$}}\;}

\def\Tilde{\widetilde}

\def\argmin{\mathop{\rm argmin}\nolimits}
\def\ra{\rangle}
\def\la{\langle}

\def\B{\Bbb B}
\def\h{\hfill\Box}
\def\R{\Bbb R}
\def\N{\Bbb N}
\def\ox{\bar{x}}
\def\ow{\bar{w}}
\def\oy{\bar{y}}
\def\oz{\bar{z}}
\def\ov{\bar{v}}
\def\ou{\bar{u}}

\def\h{\hfill\square}

\def\O{\Omega}
\def\ph{\varphi}

\newcounter{lk}

\begin{document}

\begin{center}
\vspace*{0.3in} {\bf A D.C. ALGORITHM VIA CONVEX ANALYSIS APPROACH FOR SOLVING A LOCATION PROBLEM INVOLVING SETS\footnote{The research of Nguyen Mau Nam was partially supported by the USA National Science Foundation under grant DMS-1411817 and
the Simons Foundation under grant \#208785. The research of other authors was supported by the National Foundation for Science and Technology Development (NAFOSTED), Vietnam. The authors are thankful to Prof. Nguyen Nang Tam and Mr. Nguyen Ngoc Chien for useful discusions on the subject.}}\\[2ex]
Nguyen Thai An\footnote{Thua Thien Hue College of Education, 123 Nguyen Hue, Hue City, Vietnam (thaian2784@gmail.com).}, Nguyen Mau Nam\footnote{Fariborz Maseeh Department of
Mathematics and Statistics, Portland State University, PO Box 751, Portland, OR 97207, United States (mau.nam.nguyen@pdx.edu).},
Nguyen Dong Yen\footnote{Institute
of Mathematics, Vietnam Academy of Science and Technology, 18 Hoang
Quoc Viet, Hanoi 10307, Vietnam (ndyen@math.ac.vn).}
\end{center}
{\small \textbf{Abstract.} We study a location problem that involves a weighted sum of distances to closed convex sets. As several of the weights might be negative, traditional solution methods of convex optimization are not applicable. After obtaining some existence theorems, we introduce a simple, but effective, algorithm for solving the problem. Our method is based on the Pham Dinh - Le Thi algorithm for d.c. programming and a generalized version of the Weiszfeld algorithm, which works well for convex location problems.  }

\medskip
\vspace*{0,05in} \noindent {\bf Key words.} d.c. programming, nonconvex location problem, Pham Dinh - Le Thi algorithm, difference of convex functions, Weiszfeld algorithm.

\newtheorem{Theorem}{Theorem}[section]
\newtheorem{Proposition}[Theorem]{Proposition}
\newtheorem{Remark}[Theorem]{Remark}
\newtheorem{Lemma}[Theorem]{Lemma}
\newtheorem{Corollary}[Theorem]{Corollary}
\newtheorem{Definition}[Theorem]{Definition}
\newtheorem{Example}[Theorem]{Example}
\renewcommand{\theequation}{\thesection.\arabic{equation}}
\normalsize

\section{Introduction and Problem Formulation}

The classical Fermat-Torricelli problem requires to find a point that minimizes the sum of the distances from a moving point to three given points in the plane. The first numerical algorithm for solving a generalized version of the problem that involves a finite number of points was introduced by Weiszfeld \cite{w}. Unfortunately, as shown by Kuhn \cite{k}, the algorithm may fail to converge in general. The assumptions guaranteeing the convergence of the Weiszfeld algorithm along with a proof of the convergence theorem were given in \cite{k}.

\medskip
Many generalized versions of the Fermat-Torricelli and several new algorithms have been introduced to solve generalized Fermat-Torricelli problems as well as to improve the Weiszfeld algorithm; see, e.g., \cite{b1,e,n2,mns,mv,ul,vz}. The Fermat-Torricelli problem has also been revisited several times from different viewpoints; see, e.g., \cite{b,ck,d,p,wf}.

\medskip
Motivated by applications to  more complex location problems in which the sizes of the locations are not negligible, generalized models of the Fermat-Torricelli problem involving sets have been introduced and studied intensively; see \cite{n2,mns,nh} and the references therein. To the best of our knowledge, no numerical algorithm has been developed to solve nonconvex generalized Fermat-Torricelli problems with weighted sums of distances to closed convex sets, so far. However, it is worthy to stress that the special case where all the closed convex sets reduce to single points has been treated in \cite{Tuy92}.

\medskip
In the Euclidean space $\R^n$ with the norm denoted by $\|\cdot\|$, the \emph{distance function} to a nonempty closed convex set $Q\subset \R^n$ is defined by
\begin{equation}\label{ds}
d(x;Q):=\inf\{\|x-w\|\; |\; w\in Q\}.
\end{equation}
It is well known (see e.g., \cite[the lemma on p. 53 and Proposition 2.4.1]{Cl90}) that $d(\cdot; Q)$ is a convex function with the property
$$|d(x; Q)-d(y; Q)| \leq \|x-y\|\quad ( x, y \in \R^n).$$
The unique point $w\in Q$ satisfying $\|x-w\|=d(x;Q)$ is called the \emph{Euclidean projection} of $x$ to $Q$, and is denoted by $P(x; Q)$.

\medskip
Given two finite collections $\{\O_i \; |\; i=1,\ldots, p\}$ and $\{\Theta_j \;|\;  j=1, \ldots, q\}$ of nonempty closed convex sets  in $\R^n$, we consider the constrained optimization problem
\begin{equation}\label{mainproblem}
\min\left\{f(x):=\sum_{i=1}^p\alpha_id(x; \O_i)-\sum_{j=1}^q\beta_j d(x; \Theta_j)\;|\; x\in S \right\},
\end{equation}
where $S$ is a nonempty closed convex set and the weights $\alpha_i$ and $\beta_j$ are all
positive real numbers. The convex functions $\sum_{i=1}^p\alpha_id(x; \O_i)$ and $\sum_{j=1}^q\beta_j d(x; \Theta_j)$ are two d.c. components of the d.c. function $f(x)$, where ``d.c.'' stands for ``{\bf d}ifference of {\bf c}onvex functions''. Thus, (\ref{mainproblem}) belongs to the class of nonsmooth d.c. programming problems.

\medskip
Our goal in this paper is to study (\ref{mainproblem}) from both theoretical and numerical aspects. Since the problem is nonconvex in general, traditional solution methods of convex optimization (see, e.g., \cite{r1}) are not applicable.

\medskip
First, we obtain some specific existence theorems for~\eqref{mainproblem}.  Second, to solve the problem numericaly, we introduce a simple, but effective algorithm. Our method is based on the {\em Pham Dinh - Le Thi algorithm} for d.c. programming \cite{TA1,TA2} and a generalized version of the Weiszfeld algorithm \cite{bmn}, which works well for convex location problems.

\medskip
The rest of the paper is organized as follows. Section 2 recalls several concepts and results which are used in the sequel. Solution existence theorems for~\eqref{mainproblem} in the general case are obtained in Section 3. Solution existence as well as a solution set representation, properties of the solution set and its containers (the local solution set, the stationary point set, the critical point set) of~\eqref{mainproblem} in the case $p=q=1$ are discussed in Section 4. Using the Pham Dinh - Le Thi algorithm for d.c. programming and a generalized version of the Weiszfeld algorithm, in Section 5 we develop an algorithm for solving (\ref{mainproblem}) numerically.

\section{Preliminaries}
\setcounter{equation}{0}
In this section we present some basic concepts and results used throughout the paper. The readers are referred to standard books on convex analysis such as \cite{bmn,r} for more details.

\medskip
Let $X=\R^n$ and let $Y=\R^n $ be the dual space of $X$. Both $X$ and $Y$ are equipped with the Euclidean norm in our setting. For a convex function $\ph: X \to (-\infty, \infty]$, a \emph{subgradient} of $\ph$ at $\ox\in \mbox{\rm dom}\,\ph:=\{x\in X\; |\; \ph(x)<\infty\}$ is an element $v\in Y$ such that
\begin{equation*}
\la v, x-\ox\ra\leq \ph(x)-\ph(\ox)\quad \forall \; x\in X.
\end{equation*}
The collection of all subgradients of $\ph$ at $\ox$ is called the \emph{subdifferential} of the function at this point and is denoted by $\partial \ph(\ox)$.

\medskip
According to \cite[Proposition 3.39]{bmn}, the subdifferential of the distance function \eqref{ds} at $\ox$ can be computed by the formula
\begin{equation}\label{sub_dis}
   \partial d(\ox; Q) = \left\{
     \begin{array}{lr}
       N(\ox; Q)\cap \B & \mbox{\rm if }\;\ox \in Q\\
       \left\{\dfrac{\ox- P(\ox; Q)}{d(\ox; Q)}\right\} & \mbox{\rm if }\; \ox \notin Q,
     \end{array}
   \right.
\end{equation}
where $N(\ox; Q):=\{v\in Y\; |\; \la v, x-\ox\ra\leq 0\; \mbox{\rm for all }x\in Q\}$ is the normal cone of $Q$ at $\ox$ and $\B$ stands for the closed unit ball of $\R^n$.

\medskip
The \emph{Fenchel conjugate} of a convex function $\ph: X \to (-\infty, +\infty]$ is defined by
\begin{equation*}
\ph^*(v):=\sup\{\la v,x\ra-\ph(x)\; |\; x\in \R^n\}.
\end{equation*}
By \cite[Proposition 3, p. 174]{IT79}, if $\varphi$ is proper, i.e., $\mbox{\rm dom}\,\ph\neq \emptyset$, and $\varphi$ is lower semicontinuous, then $\ph^*: Y\to (-\infty, +\infty]$ is also a proper, lower semicontinuous  convex function.

\begin{Proposition}\label{conjugate}{\rm (Properties of the Fenchel conjugates; see e.g., \cite{TA1})} Let $\ph: \R^n\to (-\infty, +\infty]$ be a convex function. \\[1ex]
{\rm (i)} Given any $x\in \mbox{\rm dom}\,\ph$, one has that $y\in \partial \ph(x)$ if and only if
\begin{equation*}
\ph(x)+\ph^*(y)=\la x, y\ra.
\end{equation*}
{\rm (ii)} If $\ph$ is proper and lower semicontinuous, then for any $x\in \mbox{\rm dom}\,\ph$ one has that $y\in \partial \ph(x)$ if and only if $x\in \partial \ph^*(y)$.\\
{\rm (iii)} If $\ph$ is proper and lower semicontinuous, then $(\ph^{*})^{*}=\ph$.
\end{Proposition}

Let $g: X\to (-\infty, +\infty]$ and $h: X \to \R$ be convex functions. Throughout the forthcoming we assume that both $g$ and $h$ are proper and lower semicontinuous. Consider the following
d.c. programming problem
\begin{equation}\label{d.c.}
\min\{f(x):=g(x)-h(x) \;|\; x\in X\}.
\end{equation}

\begin{Proposition}\label{DCO}{\rm (First-order necessary optimality condition; see e.g., \cite{TA1})}
If $\ox\in \mbox{\rm dom}\,f$ is a local minimizer of~\eqref{d.c.}, then
\begin{equation} \label{stationary}
\partial h(\ox)\subset\partial g(\ox).
\end{equation}
\end{Proposition}

Any point satisfying condition~\eqref{stationary} is called a {\em stationary point} of~\eqref{d.c.}. One says that $\ox$ a {\em  critical point} of~\eqref{d.c.} if $\partial g(\ox)\cap \partial h(\ox)\neq \emptyset$. It is obvious that every stationary point $\ox$ with $\partial h(\ox) \neq \emptyset$ is a critical point. But the converse is not true in general; see e.g., Example~\ref{4.9}.

\medskip
The \emph{Toland dual} of (\ref{d.c.}) is the problem
\begin{equation}\label{DCD}
\min\{h^*(y)-g^*(y)\;|\; y\in Y\}.
\end{equation}

Relationship between (\ref{d.c.}) and (\ref{DCD}) is described in the next proposition.
\begin{Proposition}\label{toland} {\rm (Toland's duality theorem; see e.g., \cite{TA1})} Under the assumptions made on $g$ and $h$, one has
\begin{equation*}
\inf\{g(x)-h(x)\; |\; x\in X\}=\inf\{h^*(y)-g^*(y)\; |\; y\in Y\},
\end{equation*}
i.e., the optimal values of {\rm (\ref{d.c.})} and {\rm (\ref{DCD})} coincide.
\end{Proposition}

Based on Toland's duality theorem and the results recalled in Proposition~\ref{conjugate}, Pham Dinh  and Le Thi \cite{TA1,TA2} introduced a solution method for (\ref{d.c.}) called the \emph{DCA} (d.c. algorithm). The main idea of DCA is to construct two vector sequences $\{x_k\}$ and $\{y_k\}$ such that the real sequences $g(x_k)-h(x_k)$ and $h^*(y_k)-g^*(y_k)$ are both monotone decreasing, and every cluster point $\ox$ of $\{x_k\}$ is a critical point of problem (\ref{d.c.}). Similarly, every cluster point $\oy$ of $\{y_k\}$ is a critical point of  (\ref{DCD}), i.e., $\partial g^*(\oy)\cap \partial h^*(\oy)\neq \emptyset.$

\medskip
The DCA is summarized as follows:

\smallskip
{\bf Step 1.} Choose $x_0\in \mbox{\rm dom}\,g$.

\smallskip
{\bf Step 2.} For $k\geq 0$, use $x_k$ to find $y_k\in \partial h(x_k)$. Then, use $y_k$ to find $x_{k+1}\in \partial g^*(y_k)$.

\smallskip
{\bf Step 3.} Increase $k$ by 1 and go back to {\bf Step 2}.

\medskip
In \cite{TA1}, it has been shown that the inclusion $y_k\in \partial h(x_k)$ is equivalent to the requirement that $y_k$ is a solution of the convex optimization problem
\begin{equation}\label{S1}
\min\{h^*(y)-\la x_k, y\ra\;|\; y\in Y\}.
\end{equation}
Similarly, the inclusion $x_{k+1}\in \partial g^*(y_k)$ is equivalent to the requirement that $x_{k+1}$ is a solution of the convex optimization problem
\begin{equation}\label{S2}
\min\{g(x)-\la y_k, x\ra\;|\;  x\in X\}.
\end{equation}

To solve~\eqref{S1} and~\eqref{S2}, one  has to compute the conjugate function $h^*$ and employ an appropriate solution method in convex programming. The readers are referred to \cite{TA1,TA2} for more details.

\section{Solution Existence in the General Case}
\setcounter{equation}{0}
Consider problem~\eqref{mainproblem} and put $I=\{1, \ldots, p\}, J=\{1, \ldots, q\}.$
\begin{Theorem}\label{existence}{\rm (Sufficient conditions for the solution existence)}
The problem~\eqref{mainproblem} has a solution if at least one of the following condition is satisfied:\\[1ex]
{\rm (i)} $S$ is bounded;\\
{\rm (ii)} $\sum_{i\in I}\alpha_i > \sum_{j\in J} \beta_j$, and all the sets $\O_i$, $i\in I$, are bounded.
\end{Theorem}
{\bf Proof.} (i)  The distance functions $d(\cdot; \O_i)$ and $d(\cdot; \Theta_j)$, for $i\in I$ and $j\in J$, are all Lipschitz on $\R^n$ with the Lipschitz constant 1 (see, e.g., \cite[Proposition 2.4.1]{Cl90}). Hence, $f$ is a continuous function and~\eqref{mainproblem} has a solution by the
Weierstrass theorem.

(ii) Suppose that $\sum_{i\in I}\alpha_i > \sum_{j\in J} \beta_j$ and all the sets $\O_i$ are bounded. Let $r>0$ be such that $\bigcup_{i \in I} \O_i \subset \B(0; r),$ where $\B(0; r)$ is the closed ball of radius $r$ centered at $0$. For any $i\in I$ and $x\in \R^n$, it is clear that
$$d(x; \O_i) \geq d(x; \B(0; r ))=\max\{\|x\|-r, 0\}\geq \|x\|-r.$$
For each $j\in J$, select an element $v_j\in \Theta_j$. Then, for any $x\in \R^n$,
$$d(x; \Theta_j) \leq \|x-v_j\| \quad (\forall j\in J).$$
It follows that
\begin{align*}
f(x)&=\sum_{i \in I} \alpha_i d(x; \O_i) - \sum_{j \in J} \beta_j d(x; \Theta_j)\\
&\geq \sum_{i \in I} \alpha_i (\|x\| -r) - \sum_{j \in J} \beta_j \|x-v_j\|\\
&\geq \sum_{i \in I} \alpha_i (\|x\| -r) - \sum_{j \in J} \beta_j(\|x\|+\|v_j\|)\\
&= \left (\sum_{i \in I}\alpha_i - \sum_{j \in J} \beta_j\right) \|x\|
-\left(r\sum_{i \in I}\alpha_i + \sum_{j \in J} \beta_j\|v_j\|\right).
\end{align*}
Since the last expression tends to $+\infty$ as $\|x\|\to \infty$, the objective function is coercive.
The conclusion now follows from the Weierstrass theorem and the continuity of $f$.  $\h$

\medskip
By Theorem~\ref{existence}, if $\sum_{i \in I} \alpha_i  >  \sum_{j \in J} \beta_j$ then \eqref{mainproblem} has  a solution, provided  that all the sets $\O_i$ are bounded. We are going to study the solution existence of~\eqref{mainproblem} in the cases where $\sum_{i \in I} \alpha_i  <  \sum_{j \in J} \beta_j$ or $\sum_{i \in I} \alpha_i = \sum_{j \in J} \beta_j$. In the first case, a solution exists if $S$ is bounded. Hence, just the case $S$ is unbounded needs to be considered.

\begin{Proposition}\label{existence1a} If $\sum_{i \in I} \alpha_i  <  \sum_{j \in J} \beta_j$, $S$ is unbounded, and all the sets $\Theta_j$, $j\in J$, are bounded, then $\inf\{f(x)\;|\; x\in S\}=-\infty$; so \eqref{mainproblem} has no solution.
\end{Proposition}
{\bf Proof.} Under the assumptions made, there exists $R>0$ such that $\bigcup_{j \in J} \Theta_j \subset \B(0; R).$ Given any $j\in J$ and $x\in \R^n$, we have
$$d(x; \Theta_j) \geq d(x; \B(0; R )) \geq \|x\|-R.$$
Let the elements $u_i\in\O_i$, $i\in I$, be chosen arbitrarily. Then, for every $x\in \R^n$,
\begin{align*}
f(x)&\leq \sum_{i \in I} \alpha_i \|x-u_i\| - \sum_{j \in J} \beta_j (\|x\| - R)\\
&\leq \sum_{i \in I} \alpha_i (\|x\|+\|u_i\|) - \sum_{j \in J} \beta_j (\|x\| - R)\\
&= \left (\sum_{i \in I}\alpha_i - \sum_{j \in J} \beta_j \right) \|x\|  +\left(\sum_{i \in I}\alpha_i\|u_i\| + R\sum_{j \in J} \beta_j\right).
\end{align*}
As $\sum_{i\in I}\alpha_i  < \sum_{j\in J}\beta_j$, we can assert that $\lim \limits_{\|x\|\to \infty} f(x)=-\infty.$ Hence the unboundedness of $S$ yields $\inf\{f(x)\;|\; x\in S\}=-\infty$. $\h$

\begin{Proposition}\label{existence1b}
If $\sum_{i \in I} \alpha_i  =  \sum_{j \in J} \beta_j$, and all of the sets $\O_i$,  $i\in I$ and $\Theta_j$, $j\in J$, are bounded, then there exists $\gamma >0$ such that $|f(x)|\leq \gamma$ for all $x\in \R^n$.
\end{Proposition}
{\bf Proof.} Let $r>0$ (resp., $R>0$) be chosen as in the proof of Theorem~\ref{existence} (resp., of Proposition~\ref{existence1a}). For every $j\in J$ (resp., $i\in I$), fix some $v_j\in \Theta_j$ (resp., $u_i\in \O_i$). Since $\sum_{i\in I}\alpha_i  =  \sum_{j\in J}\beta_j$, for any $x\in \R^n$, the estimates already obtained in the proofs of Theorem~\ref{existence} and
Proposition~\ref{existence1a} imply that
$$-\left(r\sum_{i \in I}\alpha_i + \sum_{j \in J} \beta_j\|v_j\|\right) \leq f(x) \leq \left(\sum_{i \in I}\alpha_i\|u_i\| + R\sum_{j \in J} \beta_j\right).$$
Setting
$$\gamma=\min\left \{r\sum_{i \in I}\alpha_i + \sum_{j \in J} \beta_j\|v_j\|, \sum_{i \in I}\alpha_i\|u_i\| + R\sum_{j \in J} \beta_j\right\},$$
we obtain the desired conclusion. $\h$

\medskip
The examples given below show that the conclusions of Theorem~\ref{existence}(ii) and Proposition~\ref{existence1a}, respectively, may not hold without the boundedness of $\O_i, i\in I$, and $\Theta_j, j\in J$.

\begin{Example}{\rm Let $n=1$, $S=\R$, $I=\{1\}$, $J=\{1\}$, $\O_1=(-\infty, 0]$, $\Theta_1=\{1\}$, $\alpha_1=2$, and $ \beta_1=1$. Since
$$f(x) = 2d(x; \O_1) - d(x; \Theta_1) \to -\infty \; \mbox{ as }  x \to - \infty,$$
problem \eqref{mainproblem} does not possess any solution.}
\end{Example}

\begin{Example}{\rm
For $n=2$, $S=\R^2$, $I=\{1\}$, $J=\{1, 2\}$, $\O_1=\R\times \{0\}$, $\Theta_1=\R\times (-\infty, -1]$, $\Theta_2=\R\times [1, \infty)$, $\alpha_1=1$, and $\beta_1=\beta_2=1$, problem~\eqref{mainproblem} becomes
$$\min\{f(x)=d(x; \O_1) -d(x; \Theta_1)-d(x; \Theta_2)\; |\; x\in \R^2\}.$$
It is not difficult to show that the solution set is $\O_1$ and the optimal value is $-2$.}
\end{Example}

If the negative weights are absent in~\eqref{mainproblem} then, as shown in \cite[Proposition 3.1]{mns}, to ensure the solution existence we only need to assume that one of the target sets $\O_i$ is bounded.

\medskip
If negative weights are present in~\eqref{mainproblem}, then the condition $\sum_{i\in I}\alpha_i > \sum_{j\in J}\beta_j$ together with the boundedness of that one of the sets $\O_i$ may be not  enough for the solution existence.
\begin{Example}{\rm
Let $n=1$, $S=\R$, $I=\{1, 2\}$, $J=\{1, 2\}$, $\O_1=\{-1\}, \O_2=[2, \infty)$, $\Theta_1=\{0\}$, $\Theta_2=\{1\}$, $\alpha_1=1$, $\alpha_2=2$, and $\beta_1=\beta_2=1$. Since for every $x\in \O_2$ we have
$$f(x)=|1+x| - |x| - |x-1|=-x +2,$$
the optimal value of~\eqref{mainproblem} is $-\infty$.}
\end{Example}

If the equality
\begin{equation}
\sum_{i\in I} \alpha_i = \sum_{j\in J} \beta_j,
\label{equalweights}
\end{equation}
holds, then the solution set of~\eqref{mainproblem} may be nonempty or empty as well. We now provide a sufficient condition for the solution existence under the assumption \eqref{equalweights}.

\begin{Proposition}\label{aux_prop}
Any solution of the problem
\begin{equation}\label{aux_mainproblem}
\max\{h(x):=\sum_{j\in J}\beta_jd(x; \Theta_j)\;|\;x\in \O_1\}
\end{equation}
is a solution of~\eqref{mainproblem} in the case where $\O_1\subset S$, $I=\{1\}$, and $\alpha_1=\sum_{j\in J} \beta_j$. Thus, in that case, if $\O_1$ is bounded then \eqref{mainproblem} has a solution.
\end{Proposition}
{\bf Proof.} Suppose $\ox\in \O_1$ is a solution of  \eqref{aux_mainproblem}. For any $x\in S$, setting $u=P(x;\O_1)$ and $v_j=P(u; \Theta_j)$ for every $j\in J$, we have
\begin{align*}
f(x)&=\alpha_1 d(x; \O_1) -\sum_{j\in J}\beta_j d(x; \Theta_j) \\
&\geq \alpha_1\|x-u\|  -\sum_{j\in J}\beta_j \|x-v_j\| \\
&= \sum_{j\in J}\beta_j\|x-u\|  -\sum_{j\in J}\beta_j \|x-v_j\| \\
&= \sum_{j\in J}\beta_j (\|x-u\|-\|x-v_j\| )\\
&\geq -\sum_{j\in J}\beta_j\|u-v_j\| =-\sum_{j\in J}\beta_jd(u; \Theta_j)=-h(u)\\
& \geq - h(\ox)=f(\ox).
\end{align*}
This shows that $\ox$ is a solution of~\eqref{mainproblem}. $\h$

\medskip

Sufficient conditions forcing the solution set of~\eqref{mainproblem} to be a subset of one of the sets $\O_i$ are given in the next proposition, which is an extension of the ``majority theorem'' in \cite{Wg} and \cite[Proposition 3]{Tuy92}.

\begin{Proposition}\label{majority}
Consider problem~\eqref{mainproblem} where $\O_{i_0} \subset S$ for some $i_0\in I$, and
$$\alpha_{i_0} > \sum_{i\in I\setminus\{i_0\}}\alpha_i + \sum_{j\in J}\beta_j.$$
Then any solution of~\eqref{mainproblem} must belong to $\O_{i_0}$.
\end{Proposition}
{\bf Proof.} Fix any $x\in S\setminus \O_{i_0}$. Let $u_i:=P(x; \O_i)$ for $i\in I$, and
$v_j:=P(u_{i_0}; \Theta_j)$ for $j\in J$. We have
\begin{align*}
f(x)&=\alpha_{i_0} d(x; \O_{i_0}) + \sum_{i\in I\setminus\{i_0\}} \alpha_i d(x; \O_i) -\sum_{j\in J}\beta_jd(x; \Theta_j)\\
&> \left(\sum_{i\in I\setminus \{i_0\}}\alpha_i + \sum_{j\in J}\beta_j \right) \|x-u_{i_0}\| +  \sum_{i\in I\setminus\{i_0\}} \alpha_i \|x-u_i\| -\sum_{j\in J}\beta_j\|x-v_j\|\\
&= \sum_{i\in I\setminus\{i_0\}}\alpha_i (\|x-u_{i_0}\| + \| x-u_i\| ) + \sum_{j\in J}\beta_j (\|x-u_{i_0}\| - \|x-v_j\|)\\
&\geq \sum_{i\in I\setminus\{i_0\}}\alpha_i \|u_{i_0} - u_i\|  - \sum_{j\in J}\beta_j \|u_{i_0}-v_j\|\\
&\geq \sum_{i\in I\setminus\{i_0\}}\alpha_i d(u_{i_0}; \O_i)  - \sum_{j\in J}\beta_j d(u_{i_0}; \Theta_j)=f(u_{i_0}).
\end{align*}
Thus, no $x\in S\setminus \O_{i_0}$ can be a solution of~\eqref{mainproblem}. $\h$

\medskip
To show that \eqref{mainproblem} can have an empty solution set under condition~\eqref{equalweights}, let us consider a special case where $S=\R^n$, $\O_i=\{a_i\}$, $\Theta_j=\{b_j\}$ with $a_i$ and $b_j$, $i\in I$ and $j\in J$, being some given points. Problem~\eqref{mainproblem} now becomes
\begin{equation}\label{prob_points}
\min\left\{f(x)=\sum_{i\in I} \alpha_i \|x-a_i\|-\sum_{j\in J}\beta_j\|x-b_j\|\;|\; x\in \R^n\right\}.
\end{equation}

\begin{Lemma}\label{opt_val_inf} Let $f(x)$ be given as in~\eqref{prob_points} and let
$w=\sum_{i\in I}\alpha_i a_i - \sum_{j\in J}\beta_j b_j $.
If $w=0$ then
\begin{equation} \label{lim0}
\lim \limits_{\|x\|\to \infty} f(x) = 0.
\end{equation}
If $w\neq 0$ then
\begin{equation} \label{liminf}
\liminf_{\|x\|\to \infty} f(x) = -\|w\|.
\end{equation}
\end{Lemma}
{\bf Proof.} By Proposition~\ref{existence1b} we can find $\gamma>0$ such  that $|f(x)| \leq \gamma$ for all $x\in \R^n$. Since
\begin{equation*}
\lim_{t\to 0}\dfrac{\sqrt{1-t}-(1-\frac{1}{2}t)}{t^2}=-1/8,
\end{equation*}
there exists $\delta\in (0, 1)$ such that
$$\left|\dfrac{\sqrt{1-t}-(1-\frac{1}{2}t)}{t^2}+\dfrac{1}{8}\right|<1$$
for any $t$ satisfying $0<|t|<\delta$. This implies that, for any $t\in(-\delta, \delta)$,
\begin{equation}\label{est1}
(1-\frac{1}{2}t)-\frac{9}{8}t^2 \leq \sqrt{1-t}  \leq (1-\frac{1}{2}t)+\frac{7}{8}t^2.
\end{equation}
It is clear that
\begin{equation}
\|x-a_i\|=\left(\|x\|^2-2\la x, a_i\ra + \|a_i\|^2\right)^{1/2} =\|x\|\left(1 - \frac{2\la x, a_i\ra}{\|x\|^2}+ \frac{\|a_i\|^2}{\|x\|^2}\right)^{1/2}.
\label{est2}
\end{equation}
Since $|\la x, a_i\ra| \leq \|x\|\|a_i\|$, we have
$$\lim \limits_{\|x\|\to \infty}\left( \frac{2\la x, a_i\ra}{\|x\|^2} - \frac{\|a_i\|^2}{\|x\|^2}\right)=0.$$
Setting $t=\frac{2\la x, a_i\ra}{\|x\|^2} - \frac{\|a_i\|^2}{\|x\|^2}$ and applying~\eqref{est1}, we can easily show that there exists $C>0$ such that
$$\left(1-\frac{\la x, a_i\ra}{\|x\|^2}\right)- \frac{C}{\|x\|^2} \leq \left(1 - \frac{2\la x, a_i\ra}{\|x\|^2}+ \frac{\|a_i\|^2}{\|x\|^2}\right)^{1/2}\leq \left(1-\frac{\la x, a_i\ra}{\|x\|^2}\right)+ \frac{C}{\|x\|^2} $$
for all $x$ with $\|x\|$ being large enough. Combining this with~\eqref{est2} we have
\begin{align}
\|x\|\left(1-\frac{\la x, a_i\ra}{\|x\|^2}\right)- \frac{C}{\|x\|}\leq \|x-a_i\|&=\left(\|x\|^2-2\la x, a_i\ra + \|a_i\|^2\right)^{1/2} \notag \\
&=\|x\|\left(1 - \frac{2\la x, a_i\ra}{\|x\|^2}+ \frac{\|a_i\|^2}{\|x\|^2}\right)^{1/2} \notag \\
&\leq \|x\|\left(1-\frac{\la x, a_i\ra}{\|x\|^2}\right)+ \frac{C}{\|x\|}, \label{estimate}
\end{align}
provided that the norm of $x$ is large enough. Hence, for such a vector $x$, for any $i\in I$ and $j\in J$ we have
$$\|x\|\left(\alpha_i - \frac{\la x, \alpha_i a_i\ra}{\|x\|^2}\right) - \frac{\alpha_iC}{\|x\|} \leq \alpha_i \|x-a_i\| \leq \|x\|\left(\alpha_i - \frac{\la x, \alpha_i a_i\ra}{\|x\|^2}\right) + \frac{\alpha_iC}{\|x\|},$$
and
$$ \|x\|\left(-\beta_j + \frac{\la x, \beta_j b_j\ra}{\|x\|^2}\right) - \frac{\beta_jC}{\|x\|} \leq  -\beta_j \|x-b_j\| \leq \|x\|\left(-\beta_j + \frac{\la x, \beta_j b_j\ra}{\|x\|^2}\right) + \frac{\beta_jC}{\|x\|}. $$
From these inequalities and the condition~\eqref{equalweights} it follows that
\begin{equation}\label{doubleinequality}
 -\bigg \la \frac{x}{\|x\|}, \sum_{i\in I}\alpha_ia_i - \sum_{j\in J}\beta_j b_j\bigg \ra -\dfrac{C_1}{\|x\|} \leq f(x) \leq -\bigg \la \frac{x}{\|x\|}, \sum_{i\in I}\alpha_ia_i - \sum_{j\in J}\beta_j b_j\bigg \ra +\dfrac{C_1}{\|x\|},
\end{equation}
where $C_1:= C\left(\sum_{i\in I}\alpha_i +\sum_{j\in J}\beta_j\right)$. Therefore, if $w=0$ then~\eqref{doubleinequality} yields~\eqref{lim0}.

Now, suppose that $w\neq 0$. By the Cauchy-Schwarz inequality we have
$$-\bigg\|  \sum_{i\in I}\alpha_ia_i - \sum_{j\in J}\beta_j b_j \bigg\| \leq -\bigg \la \frac{x}{\|x\|}, \sum_{i\in I}\alpha_ia_i - \sum_{j\in J}\beta_j b_j\bigg \ra,$$
so~\eqref{doubleinequality} implies that $$\liminf_{\|x\|\to \infty} f(x) \geq -\|w\|.$$
To obtain the equality~\eqref{liminf}, it suffices to choose $x_k=kw$ for $k\in \N$ and observe from \eqref{doubleinequality} that $\lim_{k\to \infty} f(x_k) = -\|w\|$.

The proof of the lemma is complete. $\h$

\medskip
In the next proposition, we consider problem~\eqref{prob_points} in the case where $|J|=1$.
\begin{Proposition}
Let $I=\{1, \ldots, p\},  p \geq2$, and let $b\in \R^n$. If $\beta=\sum_{i\in I}\alpha_i$ and the vectors $\{a_i -b\}$ for $i\in I$ are linearly independent, then the problem
\begin{equation}\label{prob_point_no_sol}
\min\left\{f(x)=\sum_{i\in I} \alpha_i \|x-a_i\| - \beta\|x-b\|\;|\; x\in \R^n\right\},
\end{equation}
has no solution.
\end{Proposition}
{\bf Proof. }Put $\gamma=-\|\sum_{i\in I} \alpha_ia_i -\beta b\|.$ For every $x\in \R^n$, by the triangle inequality we have
\begin{align}
f(x) - \gamma + \beta\|x-b\|&= \sum_{i\in I}\alpha_i \|x-a_i\| + \|\sum_{i\in I} \alpha_ia_i -\beta b\| \notag \\
&\geq \big \| \sum_{i\in I}\alpha_i(x-a_i) + \sum_{i\in I}\alpha_ia_i -\beta b \big\| \notag\\
&= \big \|\sum_{i\in I}\alpha_i x -\beta b  \|=\beta\|x-b\|. \label{estgamma}
\end{align}
Hence $f(x) \geq \gamma$ for all $x\in \R^n$. Since
$\liminf \limits_{\|x\| \to \infty} f(x) =\gamma$ by
Lemma~\ref{opt_val_inf}, it follows that  $\inf \limits_{x\in \R^n} f(x)=\gamma,$ i.e., the optimal value of~\eqref{prob_point_no_sol} equals to $\gamma$. To complete the proof, we have to show that $f(x)\neq \gamma$ for every $x\in \R^n$.

Suppose on the contrary that there exists $x\in \R^n$ with $f(x)=\gamma$. By the estimates~\eqref{estgamma}, $f(x)=\gamma$ if and only if
\begin{equation}
\sum_{i=1}^p\|u_i\| + \|v\| = \|\sum_{i=1}^pu_i +v\|,
\label{csinequl}
\end{equation}
where $u_i:=\alpha_i(x-a_i)$ for $i=1, \ldots, p$,  and $v:=\sum \limits_{i\in I} \alpha_ia_i -\beta b= \sum \limits_{k\in I}\alpha_k(a_k-b)$. Taking the squares of both sides of~\eqref{csinequl} and using the Cauchy-Schwarz inequality, we can deduce that $\|u_i\|\|v\|=\la u_i, v\ra$ for all $i=1,\ldots, p$ and $\|u_i\|\|u_j\|=\la u_i, u_j\ra$ for all $i, j$ with $i\neq j$. Since $v\neq 0$ by the assumed linear independence of the vectors $\{a_i-b\;|\; i=1, \ldots, p\}$, we can find nonnegative numbers $\lambda_1, \ldots, \lambda_p$ such that $u_i=\lambda_i v$ for all  $i=1,\ldots, p$. This means that
$$\alpha_i(x-a_i) =\lambda_i\sum_{k\in I}\alpha_k(a_k-b) \quad (\forall \; i\in I).$$
Therefore, $x=a_i + \frac{\lambda_i}{\alpha_i}\sum \limits_{k\in I}\alpha_k(a_k-b)$ for $i\in I$. Consequently, for any $i, j \in I$ with $i\neq j$, we have
$$a_i + \frac{\lambda_i}{\alpha_i}\sum_{k\in I}\alpha_k(a_k-b)=a_j + \frac{\lambda_j}{\alpha_j}\sum_{k\in I}\alpha_k(a_k-b).$$
Hence
\begin{equation*}
(a_i - b) - (a_j-b) + \left(\frac{\lambda_i}{\alpha_i}-  \frac{\lambda_j}{\alpha_j}\right)\sum_{k\in I}\alpha_k(a_k-b)=0.
\end{equation*}
This is equivalent to saying that
$$\left[1+\alpha_i\left( \frac{\lambda_i}{\alpha_i}-  \frac{\lambda_j}{\alpha_j} \right)\right](a_i - b) -   \left[1+\alpha_j\left( \frac{\lambda_j}{\alpha_j}-  \frac{\lambda_i}{\alpha_i} \right)\right]      (a_j-b) + \left(\frac{\lambda_i}{\alpha_i}-  \frac{\lambda_j}{\alpha_j}\right)\sum_{k\in I\setminus \{i,j\}}\alpha_k(a_k-b)=0.$$
Since the vectors $\{a_i-b\;|\; i=1, \ldots, p\}$ are linearly independent, we must have
$$\begin{cases}
1+\alpha_i\left( \frac{\lambda_i}{\alpha_i}-  \frac{\lambda_j}{\alpha_j} \right)&=0\\
1+\alpha_j\left( \frac{\lambda_j}{\alpha_j}-  \frac{\lambda_i}{\alpha_i} \right)&=0\\
\alpha_k\left(\frac{\lambda_j}{\alpha_j}-  \frac{\lambda_i}{\alpha_i}\right)&=0, \; \; \mbox{ for } k\in I\setminus\{i, j\}.
\end{cases}$$
From the first two equalities we obtain $(\alpha_i + \alpha_j)(\frac{\lambda_i}{\alpha_i}-\frac{\lambda_j}{\alpha_j})=0$, which implies that $\frac{\lambda_j}{\alpha_j}-  \frac{\lambda_i}{\alpha_i}=0$. Substituting this back to the equality $1+\alpha_i\left( \frac{\lambda_i}{\alpha_i}-  \frac{\lambda_j}{\alpha_j} \right)=0$, we get a contradiction. $\h$

\section{Solution Existence in a Special Case}
\setcounter{equation}{0}
Consider a special case of problem~\eqref{mainproblem} where $p=q=1$ and  $S=\R^n$;
that is,
\begin{equation}\label{problem}
\min\{f(x):=\alpha d(x; \O) - \beta d(x; \Theta)\;|\; x\in \R^n\},
\end{equation}
where $\alpha \geq \beta>0$.  Let $g(x):=\alpha d(x; \O)\; \mbox{ and }\; h(x):=\beta d(x; \Theta).$

\medskip
We are going to establish several properties of the optimal solutions to problem (\ref{problem}).
The relationship between \eqref{problem} and the problem
\begin{equation}\label{problemnew}
\max\{d(x; \Theta)\;|\; x\in \O\}
\end{equation}
will be also discussed.
\begin{Lemma}\label{pp1} There is no local solution of~\eqref{problem} that belongs to the set $\Theta\setminus \O$.
\end{Lemma}
{\bf Proof.} Suppose by contradiction that there exists $\ox \in \Theta\setminus \O$ which is a local solution of \eqref{problem}. Let $\delta>0$ be such that $f(x) \geq f(\ox)$ for all $x\in\B(\ox; \delta).$ Setting $\ow:=P(\ox; \O)$, we have $\ow \neq \ox$. For $x_t:=\ox +t(\ow -\ox)$ with $t\in (0,1)$ being sufficiently small, one has
\begin{align*}
f(x_t) &= \alpha d(x_t; \O) - \beta d(x_t; \Theta) = \alpha (1-t)\|\ox-\ow\| - \beta d(x_t; \Theta) \\
&< \alpha\|\ox-\ow\| = \alpha d(\ox; \O)=f(\ox).
\end{align*}
Since $ x_t \in \B(\ox; \delta)$ for $t>0$ small enough, we have arrived at a contradiction.$\h$
\begin{Proposition}\label{pp2} If $\alpha >\beta$ then any local solution of~\eqref{problem} must belong to $\O$.
\end{Proposition}
{\bf Proof.} If $\ox$ is a local optimal solution of  \eqref{problem} and $\ox\notin \O$, then $\ox \notin \O\cup \Theta$ by Lemma \ref{pp1}. Employing Proposition \ref{DCO} we have $\partial h(\ox) \subset \partial g(\ox).$
Since $\ox \notin \O\cup \Theta$, it follows from formula \eqref{sub_dis} and the last inclusion that
$$\beta \dfrac{\ox - P(\ox; \Theta)}{d(\ox; \Theta)} =  \alpha \dfrac{\ox - P(\ox; \O)}{d(\ox; \O)}.$$
By taking the norms of the vectors on both sides, we have $\alpha=\beta$. This contradicts the assumption $\alpha >\beta$. $\h$
\begin{Proposition}\label{pp3} If $\alpha > \beta$, then $\ox$ is a solution of~\eqref{problem} if and only if it is a solution of \eqref{problemnew}. Thus, in the case $\alpha > \beta$, the solution set of \eqref{problem} does not depend on the choice of $\alpha$ and $\beta$.
\end{Proposition}
{\bf Proof.} If $\ox$ is a solution of~\eqref{problem}, then $\ox\in \O$ by Proposition \ref{pp2}.
Hence, for every $x\in \R^n$
\begin{equation*}
\alpha d(x;\O)-\beta d(x; \Theta)\geq \alpha d(\ox; \O)-\beta d(\ox;\Theta)=-\beta d(\ox; \Theta).
\end{equation*}
This implies $-d(x; \Theta)\geq -d(\ox; \Theta)$ for all $x\in \O$; so $\ox$ is a solution of \eqref{problemnew}.

Conversely, suppose that $\ox$ is a solution of \eqref{problemnew}. Then $\ox\in \O$ and $d(x; \Theta)\leq d(\ox; \Theta)$ for all $x\in \O$. This implies $f(x)\geq f(\ox)$ for all $x\in \O$.
Given any $x\notin \O$, we put $u=P(x;\O)$. By the estimates obtained in the proof of Proposition~\ref{majority}, we have $f(u)<f(x)$. Therefore,
\begin{equation*}
f(x) > f(u) =-\beta d(u;\Theta) \geq -\beta d(\ox; \Theta)=\alpha d(\ox; \O)-\beta d(\ox; \Theta)=f(\ox).
\end{equation*}
Thus we have $f(x)\geq f(\ox)$ for all $x\in \R^n$. So $\ox$ is a solution of  \eqref{problem}. $\h$
\begin{Lemma}\label{pp4} If $\alpha=\beta$ and if \eqref{problem} has a local solution $\ox$ with $\ox\notin \O$, then there must exist $\ou\in \O$ with $f(\ou)\leq f(\ox)$.
\end{Lemma}
{\bf Proof.} Since $\alpha=\beta>0$, there is no loss of generality in assuming that $\alpha=\beta=1.$ Let $\ox$ be a local solution of \eqref{problem} with $\ox\notin \O$. By Lemma \ref{pp1}, $\ox\notin \O\cup \Theta$. It follows from Proposition \ref{DCO} and formula \eqref{sub_dis} that
\begin{equation}\label{dkcan}
\dfrac{\ox - P(\ox; \Theta)}{d(\ox; \Theta)} =  \dfrac{\ox - P(\ox; \O)}{d(\ox; \O)}.
\end{equation}
Setting $\ou=P(\ox; \O)$, $\ov=P(\ox; \Theta)$, and $ \lambda=\dfrac{d(\ox; \O)}{d(\ox; \Theta)},$ by \eqref{dkcan} we have $\ou=(1-\lambda)\ox + \lambda \ov$. If $\lambda >1$ then
\begin{equation*}
f(\ox) = d(\ox; \O) -d(\ox; \Theta) >0\geq -d(\ou; \Theta) =f(\ou).
\end{equation*}
Now, suppose that $\lambda \leq 1$. Since $\ov=P(\ox; \Theta)$, for any $v\in \Theta$ we have
\begin{align*}
\la \ou -\ov, v-\ov \ra & = \la (1-\lambda)\ox + \lambda \ov -\ov, v-\ov\ra\\
&=(1-\lambda)\la \ox -\ov, v-\ov\ra\leq 0.
\end{align*}
Therefore $P(\ou; \Theta)=\ov$, and hence
$$f(\ox)=d(\ox; \O)-d(\ox; \Theta) = \|\ox - \ou\| - \|\ox -\ov\| \geq - \|\ou-\ov\|=-d(\ou; \Theta)=f(\ou).$$
We have thus found $\ou\in \O$ with $f(\ou)\leq f(\ox)$. $\h$

\begin{Proposition}\label{pp5} Suppose $\alpha=\beta$. If $\bar{u}$ is a solution of \eqref{problemnew}, then $\bar{u}$ is a solution of \eqref{problem}. Conversely, if \eqref{problem} has a solution, then \eqref{problemnew} also has a solution {\rm(}but the solution set of ~\eqref{problem} may be larger than that of \eqref{problemnew}{\rm)}.
\end{Proposition}
{\bf Proof.} We may assume that $\alpha=\beta=1.$ The first assertion follows from Proposition~\ref{aux_prop}. To prove the second one, suppose that  $\ox$ is a solution of  \eqref{problem}. Then
 $$d(\ox;\O)-d(\ox;\Theta)=f(\ox) \leq f(x)=d(x; \O)-d(x;\Theta)$$
for all  $x\in \R^n$. If $\ox\in \O$, then this implies that $\ox$ is a solution of \eqref{problemnew}. Otherwise, by Lemma \ref{pp4} we can find $\ou\in \O$ such that $f(\ou)\leq f(\ox)$. Thus,
$f(\ou) = f(\ox) \leq f(u)$ for every $u\in \O$; hence $\ou$ is a solution of \eqref{problemnew}.  $\h$

\medskip
We now describe a relationship between the solution sets of \eqref{problem} and \eqref{problemnew}, which are denoted respectively by $S_1$ and $S_2$.
\begin{Proposition}\label{pp6} Suppose that $\O \setminus \Theta \neq \emptyset$. If $\alpha=\beta$, then
\begin{equation}
S_1= \left\{\bar{u}+\R^{+}\left(\bar{u} -P(\bar{u}; \Theta)\right)\; \big|\;  \bar{u} \in S_2\right\}.
\label{conection_two_sols}
\end{equation}
\end{Proposition}
\noindent {\bf Proof.} Again, we assume that $\alpha=\beta=1$. Fix any $\bar{x}\in S_1$. If $\bar{x}\in \O$, then $\bar{x}\in S_2$. Consider the case $\ox \notin \O$. Since $\O \setminus \Theta \neq \emptyset$,  we can select a vector $w \in \O \setminus \Theta$. Observe that
\begin{equation}
f(\ox) = d(\ox; \O) - d(\ox; \Theta) \leq f(w) = -d(w; \Theta) <0.
\label{lamda}
\end{equation}
Let $\bar{u}:=P(\bar{x}; \O)$, $\bar{v}:=P(\bar{x}; \Theta)$, and $\lambda:=\dfrac{d(\ox; \O)}{d(\ox; \Theta)}$. As shown in the proof
of Lemma \ref{pp4}, $\ov=P(\ou; \Theta)$ and $f(\bar{u})\leq f(\bar{x})$. Since $\ox\in S_1$, this yields $f(\bar{u})=f(\bar{x})$. Then we have $\bar{u}\in S_2$. Moreover, it follows from~\eqref{lamda} that $0 \leq \lambda < 1$. By \eqref{dkcan} we obtain $\bar{x}=\bar{u} + \frac{\lambda}{1-\lambda}(\bar{u} - \bar{v})$; hence $\bar{x}$ belongs to the set on the right-hand side of~\eqref{conection_two_sols}.

To complete the proof, it suffices to show that any vector of the form $\bar{x}=\bar{u} + t(\bar{u} - \bar{v})$ where $\bar{u} \in S_2$, $\bar{v}=P(\bar{u}; \Theta)$, and $t\geq 0$, belongs to $S_1$. Fix such an $\ox$. Then, for any $v\in \Theta$ we have
\begin{align*}
\la v-\bar{v}, \bar{x} -\bar{v}\ra &= \la v-\bar{v}, \bar{u} +t(\bar{u}-\bar{v}) -\bar{v}\ra\\
&=(t+1)\la v-\bar{v}, \bar{u} -\bar{v}\ra\leq 0.
\end{align*}
Thus, $\bar{v}=P(\bar{x};\Theta )$. Since $\bar{u} \in S_2$, by Proposition \ref{pp5} we have $\bar{u} \in S_1$.  Therefore, with $\hat{u}:=P(\bar{x}; \O)$,
$$\|\bar{x}-\hat{u}\| -\|\bar{x} -\bar{v}\| = f(\bar{x}) \geq f(\bar{u})=   -\|\bar{u}-\bar{v}\|.$$
Then
\begin{align*}
\|\bar{x}-\hat{u}\|  &\geq \|\bar{x} -  \bar{v}\| - \|\bar{u}-\bar{v}\| \\
&=\|\bar{u} + t(\bar{u} -\bar{v}) -  \bar{v}\| - \|\bar{u} -\bar{v}\| = t\|\bar{u} -\bar{v}\|=\|\bar{x} -\bar{u}\|.
\end{align*}
Since $\hat{u}=P(\bar{x}; \O)$, this implies that $\hat{u}= \bar{u}$. Now we have
\begin{align*}
f(\bar{x})=\|\bar{x} -\hat{u} \| -\|\bar{x}-\bar{v}\| &=\|\bar{x} -\bar{u} \| -\|\bar{x}-\bar{v}\| \\
&= t\|\bar{u} -\bar{v}\| - (t+1)\|\bar{u} -\bar{v}\|\\
&=-\|\bar{u}- \bar{v}\|=f(\bar{u}).
\end{align*}
Recalling that $\ou\in S_1$, from this we obtain $\bar{x} \in  S_1$, as desired. $\h$
\begin{Proposition}\label{pp7} Problem~\eqref{problem} has a unique solution if and only if one of the following conditions holds:\\
{\rm (i)} $\alpha>\beta$ and problem~\eqref{problemnew} has a unique solution.\\
{\rm (ii)}   $\alpha = \beta$ and $\O$ is a singleton contained in the interior of $\Theta$.
\end{Proposition}
{\bf Proof.} First, consider the case $\alpha>\beta$. Since $S_1=S_2$ by Proposition \ref{pp3},
\eqref{problem} has a unique solution if and only if \eqref{problemnew} has a unique solution.

Next, suppose that $\alpha =\beta$.

If \eqref{problem} has a unique solution then, by Proposition \ref{pp5}, problem~\eqref{problemnew} has a solution $\bar{u}$. If $\O \setminus \Theta \neq \emptyset$, then  $\bar{u}\notin \Theta$. Indeed, if $\bar{u}\in \Theta$ then for any $u\in \O$ we have $d(u; \Theta) \leq d(\bar{u}; \Theta)=0$. This yields $\O\subset \Theta$, which contradicts the assumption $\O \setminus \Theta \neq \emptyset$. Now, the property $\bar{u}\notin \Theta$ implies that $P\left(\bar{u}; \Theta\right)\neq \bar{u}$. By the representation~\eqref{conection_two_sols}, $S_1$ contains a solution ray, contradicting the solution uniqueness of \eqref{problem}. We have shown that $\O\subset \Theta$. The last inclusion gives us
$$f(x)=d(x; \O) - d(x; \Theta) \geq 0$$
for any $x\in \R^n$. The equality $f(x)=0$ holds if and only if $d(x; \O) = d(x; \Theta)$.
Since $f(x)=0$ for any $x\in \O$, we have $\O \subset S_1$. Moreover, as $S_1\cap \O\subset S_2$, we can assert that   $\O \subset S_1\cap S_2$. Therefore $\O$ must be a singleton.
If $\O=\{u\}$ and $u$ is a boundary point of $\Theta$, then $N(u; \Theta)\neq \{0\}$; see, e.g., \cite[Corollary 2.12]{bmn}. Taking any $v\in N(u; \Theta)\setminus\{0\}$, we observe that $P(u + tv; \Theta)=\{u\}=P(u + tv; \O)$ for every $t\geq 0$. Hence $f(u + tv)=0$ for any $t\geq 0$. It follows that $S_1$ contains the solution ray $u + \R^{+}v$, contrary to the assumed solution uniqueness of  \eqref{problem}. Thus $u$ must belong to the interior of $\Theta$.

Now, let $\O=\{u\}$ and $u$ is an interior point of $\Theta$. It is easy to show that $f(u)=0$ and $f(x)>0$ for every $x\in \R^n\setminus\{u\}$. Hence $u$ is a unique solution of \eqref{problem}.

The proof is complete. $\h$

\begin{Example}{\rm
If $\O\subset \Theta$, then it is clear that $S_2=\O$. Moreover, the set on the right-hand side of~\eqref{conection_two_sols} is $S_2$. But the equality $S_1=S_2$ does not hold in general. For example, in $\R^2$, let $\O$ be the closed ball centered at $(-1,0)^\top$ with radius $1$ and let $\Theta$ be the closed ball centered at $(-2, 0)^\top$ with radius 2, where the superscript
$^\top$ denotes matrix transposition. Then we have
 $$S_2=\O \subset \left[ \O \cup (\R^{+}\times\{0\})\right]=S_1.$$}
\end{Example}

Let us analyze furthermore the solution structure of~\eqref{problem}. For our convenience, we denote the {\em solution set}, the {\em local solution set}, the {\em stationary point set} and the {\em critical point set} of that problem respectively by ${\rm Sol}(P)$, ${\rm loc}(P)$, ${S}(P)$, and ${S_0}(P)$. Thus, in the preceeding notation, we have ${\rm Sol}(P)=S_1.$ Moreover, it is clear that
\begin{equation}\label{inclusions}
{\rm Sol}(P) \subset {\rm loc}(P) \subset {S}(P) \subset {S_0}(P).
\end{equation}

By constructing suitable examples, we are going to show that all the three inclusions in~\eqref{inclusions} can be strict, and

\qquad- ${\rm Sol}(P)$ can be nonconvex,

\qquad - ${\rm loc}(P)$ can be nonconvex and nonclosed,

\qquad - $S(P)$ can be nonconvex and nonclosed,

\qquad - $S_0(P)$ can be nonconvex.

\medskip
In the following examples, it  is always assumed that $n=2$ and $\alpha=\beta=1$.
\begin{Example}\label{4.9}{\rm
Let $\O=\{0\}\times [-2, 2]$ and $\Theta=\{(1, 0)^\top\}$. It is easy to see that the solution set of  problem~\eqref{problemnew} is
$$S_2=\bigg\{(0, -2)^\top, \; (0, 2)^\top \bigg\}.$$
According to~\eqref{conection_two_sols},
$${\rm Sol}(P)=\left\{(0, -2)^\top -t(1, 2)^\top\;|\; t\geq 0 \right\} \cup \left\{(0, 2)^\top +t(-1, 2)^\top\;|\; t\geq 0 \right\}.$$
By definition, $\ox\in S(P)$ if and only if $\partial d(\cdot; \Theta)(\ox)\subset \partial d(\cdot; \O)(\ox)$. From the subdifferential formula \eqref{sub_dis} we can deduce that
$$S(P)={\rm Sol}(P)\cup \bigg( (-\infty, 0] \times \{0\}\bigg) \cup \bigg( (1, +\infty)\times \{0\}\bigg).$$
If we can show that every point of the form $x=(x_1, 0)^\top \in S(P)\setminus {\rm Sol}(P)$ doesn't belong to ${\rm loc}(P)$, then combining this with ${\rm Sol}(P) \subset {\rm loc}(P) \subset {S}(P)$ we get ${\rm loc}(P)={\rm Sol}(P)$. For $\varepsilon > 0$ small enough, let $x_\varepsilon=(x_1,\varepsilon)^\top$. Since
$$f(x_\varepsilon)=d(x_\varepsilon ; \O) - d(x_\varepsilon; \Theta) = d\left(x; \O\right)- \sqrt{\epsilon^2 + d^2(x; \Theta)} < d(x; \O) - d(x; \Theta)=f(x),$$
we can assert that $x\notin {\rm loc}(P)$. To find $S_0(P)$, it suffices to observe from~\eqref{sub_dis} that any critical  point $\ox$ with $\ox\notin \O\cup \Theta$ must be a stationary point.  So we have
$$S_0(P)=S(P) \cup \Theta=S(P) \cup \{(1, 0)^\top\}.$$
}
\end{Example}

\begin{Example}\label{4.10}{\rm
Let $\O=\{0\}\times [-2, 1]$ and $\Theta=\{(1, 0)^\top\}$. Here we have
\begin{align*}
{\rm Sol}(P)&=\left\{(0, -2)^\top -t(1,2)^\top\;|\; t\geq 0 \right\},\\
{\rm loc}(P)&=\left\{(0, -2)^\top -t(1, 2)^\top\;|\; t\geq 0 \right\} \cup \left\{(0, 1)^\top +t(-1, 1)^\top\;|\; t\geq 0 \right\},  \\
S(P)&={\rm loc}(P)\cup \bigg( (-\infty, 0] \times \{0\}\bigg) \cup \bigg( (1, +\infty)\times \{0\}\bigg), \\
S_0(P)&=S(P)\cup \big\{(1, 0)^\top\big\}.
\end{align*}
In particular, ${\rm Sol}(P)\neq {\rm loc}(P)$.
}
\end{Example}

\begin{Example}\label{4.11}{\rm
Let $\O=\{0\}\times [-2, 2]$ and $\Theta=\{1\}\times [-1, 1]$. In this problem,
\begin{align*}
{\rm Sol}(P)&= \left\{(0, -2)^\top -t(1, 1)^\top\;|\; t\geq 0 \right\} \cup \left\{(0, 2)^\top +t(-1, 1)^\top\;|\; t\geq 0 \right\}, \\
{\rm loc}(P)&={\rm Sol}(P) \cup \bigg ( (-\infty, 0]\times (-1, 1)\bigg) \cup \bigg( (1, +\infty)\times (-1, 1)\bigg),\\
S(P)&={\rm loc}(P)\cup \bigg( (-\infty, 0] \times \{-1, 1\}\bigg) \cup \bigg( (1, +\infty)\times \{-1, 1\}\bigg), \\
S_0(P)&=S(P)\cup \bigg(\{1\}\times [-1, 1]\bigg).
\end{align*}
In particular, ${\rm loc}(P)$ can be nonconvex and nonclosed.
}
\end{Example}

\section{DCA and Generalized Weiszfeld Algorithm}
\setcounter{equation}{0}

To obtain  a simple algorithm for solving the generalized Fermat-Torricelli problem \eqref{mainproblem} by the DCA, we rewrite  \eqref{mainproblem}  equivalently as
\begin{equation}\label{DCprob}
\min\left\{\sum_{i\in I} \alpha_id(x;\O_i)+\delta(x; S)-\sum_{j\in J}\beta_jd(x; \Theta_j)\;|\; x\in \R^n\right\},
\end{equation}
where $\delta(\cdot; S)$ is the \emph{indicator function} of $S$ defined by $\delta(x;S)=0$ if $x\in S$, and $\delta(x;S)=\infty$ otherwise. The objective function of~\eqref{DCprob} can be put in the form $g(x)- h(x)$ with
$$g(x):=\sum_{i\in I}\alpha_id(x; \O_i) + \frac{\lambda}{2}\|x\|^2+\delta(x; S),
\quad h(x):=\sum_{j\in J}\beta_j d(x; \Theta_j) + \frac{\lambda}{2}\|x\|^2,$$
and $\lambda > 0$ being an arbitrarily chosen constant. The DCA (see the final part of Section 2) is applicable to the problem
$$\min\{g(x)-h(x)\;|\; x\in \R^n\}.$$

It follows from~\eqref{sub_dis} that if $x_k\notin \Theta_j$ then
$$\partial d(x_k; \Theta_j)=\left \{\dfrac{x_k-P(x_k;\Theta_j)}{d(x_k;\Theta_j)}\right \},$$
otherwise $0\in \partial d(x_k; \Theta_j)$. Therefore, as an element $y_k\in \partial h(x_k)$, we can choose $y_k = \sum_{j\in J} \beta_j w_j  + \lambda x_k$, where $w_j:=\dfrac{x_k-P(x_k;\Theta_j)}{d(x_k;\Theta_j)}$ if $x_k\notin \Theta_j$ and $w_j=0$ otherwise.

\medskip
With $v:=y_k$, the auxiliary problem~\eqref{S2} becomes
\begin{align*}
(P_v)\hspace{2cm} &\min\left\{\sum_{i\in I} \alpha_id(x;\O_i)+ \frac{\lambda}{2}\|x\|^2-\la v, x\ra\;|\; x\in S\right\}.
\end{align*}
We can  solve this problem by several methods of convex programming. Our solution method for $(P_v)$ is based on the Weiszfeld algorithm; see \cite{k,w}.

\medskip
For simplicity, assume that $\O_i\cap S=\emptyset$ for every $i\in I$. Fix an element $v\in \R^n$ and consider the function
\begin{equation*}
\ph_v(x):=\sum_{i\in I} \alpha_id(x;\O_i)+ \frac{\lambda}{2}\|x\|^2-\la v, x\ra,\quad x\in \R^n.
\end{equation*}
By formula~\eqref{sub_dis}, for each $x\in S$ we have
\begin{equation*}
\nabla\ph_v(x):=\sum_{i\in I} \alpha_i\dfrac{x-P(x;\O_i)}{d(x;\O_i)}+\lambda x-v.
\end{equation*}
Solving the equation $\nabla \ph_v(x)=0$ yields
\begin{equation*}
x=\dfrac{\sum \limits_{i\in I}\dfrac{\alpha_i P(x;\O_i)}{d(x;\O_i)} + v}{\sum \limits_{i\in I}\dfrac{\alpha_i}{d(x;\O_i)}+\lambda}.
\end{equation*}
Define
\begin{equation}\label{fv}
F_v(x)=\dfrac{\sum \limits_{i\in I}\dfrac{\alpha_i P(x;\O_i)}{d(x;\O_i)} +  v}{\sum \limits_{i\in I}\dfrac{\alpha_i}{d(x;\O_i)}+\lambda}, \quad x\in S.
\end{equation}
We introduce the following \emph{generalized Weiszfeld algorithm} to solve $(P_v)$:
\begin{itemize}
\item Choose $x_0\in S$.
\item Find $x_{k+1}=P(F_v(x_k); S)$ for $k\in \N$, where the mapping $F_v$ is defined in \eqref{fv}.
\end{itemize}

\begin{Theorem}\label{decreasing}\label{p2}
Consider the generalized Weiszfeld algorithm for solving $(P_v)$. If $x_{k+1}\neq x_k$, then $\varphi_v(x_{k+1}) < \varphi_v(x_{k})$.
\end{Theorem}
{\bf Proof.} Given any $x\in S$, we define
$$g(z)=\sum_{i\in I} \dfrac{\alpha_i\|z-P(x;\O_i)\|^2}{d(x;\O_i)}+\lambda\|z\|^2 -2\la v, z\ra.$$
Then
$$\nabla g(z)=2\left(\sum_{i\in I} \alpha_i\dfrac{z-P(x;\O_i)}{d(x;\O_i)} + \lambda z-v\right).$$
Since $g(z)$ is strongly convex, it has a unique global minimizer  $\oz$ on $S$ which satisfies
the inclusion
\begin{equation}
-\nabla g(\oz)\in N(\oz; S).
\label{nabla_g}
\end{equation}
Observe that
$$-\nabla g(\oz) = 2\left(\sum_{i\in I}\frac{\alpha_i}{d(x; \O_i)}+\lambda\right)(F_v(x)-\oz).$$
Hence, from~\eqref{nabla_g} it follows that
\begin{equation*}
F_v(x)-\oz\in N(\oz; S).
\end{equation*}
Therefore, $P(F_v(x); S)=\oz$.  If $\oz \neq x$, we have strictly inequality $g(\oz) < g(x)$.
Since $\|x-P(x; \O_i)\|=d(x; \O_i)$ and $\|\oz-P(x; \O_i)\| \geq \|\oz-P(\oz; \O_i)\|$ for every $i \in I$, we have
\begin{align*}
g(\oz)&=\sum_{i\in I} \dfrac{\alpha_i\|\oz-P(x;\O_i)\|^2}{d(x;\O_i)}+\lambda\|\oz\|^2 -2\la v, \oz\ra\\
&=\sum_{i\in I} \alpha_i \dfrac{\big(\|x-P(x;\O_i)\|+ \|\oz-P(x;\O_i)\| - \|x-P(x;\O_i)\|  \big)^2}{d(x;\O_i)} +\lambda\|\oz\|^2 -2\la v, \oz\ra\\
&= 2\left( \sum_{i\in I} \alpha_i\|\oz-P(x;\O_i)\|+\frac{\lambda}{2}\|\oz\|^2-\la v, \oz\ra \right)
- \left(\sum_{i\in I} \alpha_i\|x-P(x;\O_i)\|+\frac{\lambda}{2}\|x\|^2-\la v, x\ra \right) \\
&+\sum_{i\in I} \alpha_i \dfrac{\big(\|\oz-P(x;\O_i)\| - \|x-P(x;\O_i)\|\big)^2}{d(x; \O_i)}+\frac{\lambda}{2}\|x\|^2 - \la v,x\ra\\
&\geq 2\varphi_v\left(\oz\right) - \varphi_v(x) +\sum_{i\in I} \alpha_i \dfrac{\big(\|\oz-P(x;\O_i)\| - \|x-P(x;\O_i)\|\big)^2}{d(x;\O_i)} +\frac{\lambda}{2}\|x\|^2- \la v,x\ra.
\end{align*}
As  $g(\oz) < g(x)=\sum_{i\in I} \alpha_i\|x-P(x;\O_i)\|-2\la v, x\ra +\lambda\|x\|^2 = \varphi_v(x) -\la v,x\ra +\frac{\lambda}{2}\|x\|^2$, it holds that $$\varphi_v(\oz) < \varphi_v(x).$$
Substituting $x=x_k$ and $\oz=x_{k+1}=P(F_v(x_k); S)$ into the last inequality yields $\varphi_v(x_{k+1})< \varphi_v(x_k)$, provided that $x_{k+1}\neq x_k$. $\h$

\medskip
Our convergence theorem  for the generalized Weiszfeld algorithm can be formulated as follows
\begin{Theorem} The sequence $\{x_k\}$ produced by the generalized Weiszfeld algorithm converges to the unique solution of problem $(P_v)$.
\end{Theorem}
{\bf Proof.} Since $\ph_v(x)$ is strongly convex, problem $(P_v)$ has a unique solution. By the coerciveness of $\ph_v(x)$, there exists $r>0$ such that
$$\ph_v(x) > \ph_v(x_0) \; \; \mbox{ whenever }\; \|x\| >r.$$
According to Theorem~\ref{decreasing}, $\ph_v(x_k) \leq \ph_v(x_0)$. This implies
$\|x_k\| \leq r$ for all $k$, and thus the sequence $\{x_k\}$ is bounded. Let $\{x_{k_\ell}\}$ be an arbitrary subsequence of $\{x_k\}$ with a limit $\ox\in S$. The mapping $\psi(x): S\to S$ defined by $\psi(x):=P(F_v(x); S)$ and the function $\ph_v$ satisfy the assumptions of \cite[Proposition~1]{CZL12}, so $\|x_{k+1}-x_k\|\to 0$ as $k\to \infty$. The continuity of $F_v(\cdot)$ and the projection mapping implies that $\psi$ is continuous. Since
$$x_{k_\ell+1}=\psi(x_{k_\ell}),$$
one obtains $\psi(\ox)=P(F_v(\ox); S)=\{\ox\}$, and hence $\ox$ is the unique solution of $(P_v)$. Indeed, we have $F_v(\ox)-\ox\in N(\ox; S)$, which implies $-\nabla \ph_v(\ox)\in N(\ox; S)$. Therefore, $\ox$ coincides with the unique solution of $(P_v)$.
It follows that the entire sequence $\{x_k\}$ converges to $\ox$. $\h$

\medskip
Combining the DCA presented in Section 2 and the generalized Weiszfeld algorithm, we get the following algorithm for solving \eqref{mainproblem}.

\smallskip
{\bf Step 1A.} Choose $x_0\in S$, $\lambda> 0$ and $K\in \N$.

\smallskip
{\bf Step 2A.}  For $k\in \{1, \ldots, K\}$:

\begin{itemize}
\item[-] Find $y_k = \sum_{j\in J} \beta_j w_j  + \lambda x_k$, where $w_j:=\dfrac{x_k-P(x_k;\Theta_j)}{d(x_k;\Theta_j)}$ if $x_k\notin \Theta_j$ and $w_j=0$ otherwise.
\item[-] Find the unique solution $x_{k+1}= \argmin \limits_{x\in S} \varphi_{y_k}(x)$ of $(P_v)$ where $v:=y_k$ by the generalized Weiszfeld algorithm, provided that a stopping criterion and  a starting point $z_k$ are given.
\end{itemize}

\smallskip
{\bf Step 3A.}  Output: $x_{K+1}$.

\begin{Theorem}
If either condition {\rm (i)} or {\rm (ii)}  in Theorem~\ref{existence} is satisfied, then any limit point of the interative sequence $\{x_k\}$ is a critical point of the problem~\eqref{DCprob}.
\end{Theorem}
{\bf Proof.} Let $f(x):=\Tilde{g}(x)-\Tilde{h}(x)$, where $\Tilde{g}(x):=\sum_{i\in I} \alpha_i d(x; \O_i)$ and $\Tilde{h}(x):=\sum_{j\in J} \beta_j d(x; \Theta_j)$.
Recall that
$$\varphi_{y_k}(x) = \sum_{i\in I} \alpha_i d(x; \O_i) -\la y_k, x\ra +\frac{\lambda}{2}\|x\|^2.$$
On one hand, since $x_{k+1}=\argmin \limits_{x\in S} \varphi_{y_k}(x)$, we have
$$  - \nabla \varphi_{y_k}(x_{k+1}) \in N(x_{k+1}; S).$$
This means that
\begin{equation}
y_k - \lambda x_{k+1} - \nabla \Tilde{g}(x_{k+1}) \in N(x_{k+1}; S).
\label{partial_g}
\end{equation}
On the other hand, the formula for $y_k$ in {\bf Step 2A} yields
\begin{equation}
y_k - \lambda x_k \in \partial \Tilde{h}(x_k).
\label{partial_h}
\end{equation}
Hence, from the inclusion $\partial d(x_{k}; \Theta_j) \subset \B$ we can deduce
\begin{equation}
\|y_k\| \leq \lambda\|x_k\| +\sum_{j\in J}\beta_j.
\label{yk_bounded}
\end{equation}
\indent If condition {\rm (i)} in Theorem~\ref{existence} is satisfied, then by the boundedness of $S$ and by~\eqref{yk_bounded} we can assert that the sequences $\{x_k\}$ and $\{y_k\}$ are bounded.

If condition {\rm (ii)} in Theorem~\ref{existence} is satisfied then, as it has been shown in the proof of Theorem~\ref{existence}, $f(x)$ is coercive.
Since $h(x)$ is strongly convex with modulus $\lambda$, by \cite[Theorem 3.7 (i)]{TA2} we have
$$f(x_{k+1}) \leq f(x_k) - \frac{\lambda}{2}\|x_{k+1}-x_k\|^2.$$
This descent property of $\{f(x_k)\}$ in combination with the coerciveness of $f(x)$ imply that
$\{x_k\}$ is bounded.  Then $\{y_k\}$ is also bounded by~\eqref{yk_bounded}.

Consequently, there exist subsequences $\{x_{k_\ell}\}$ and $\{y_{k_\ell}\}$ of $\{x_k\}$ and $\{y_k\}$, respectively, such that $x_{k_\ell} \to \ox$ and $y_{k_\ell}\to \oy$. Since $\lim \limits_{l\to \infty} \|x_{k_\ell+1} - x_{k_\ell}\| =0$ by \cite[Theorem 3.7 (iii)]{TA2}, we have
$\lim \limits_{\ell\to \infty} \|x_{k_\ell+1} - \ox\| =0$.

Passing the inclusions~\eqref{partial_g} and~\eqref{partial_h} to limits as $k=k_\ell\to \infty$ implies
$$\oy -\lambda\ox \in \partial \Tilde{h}(\ox) \cap \left(\nabla \Tilde{g}(\ox) + N(\ox; S)\right).$$
This amounts to saying that $\oy \in \partial h(\ox) \cap \partial g(\ox).$ Thus, $\ox$ is a critical point of~\eqref{DCprob}. $\h$

\begin{Example}\label{ExPoint}{\rm The latitude/longitude coordinates in decimal format of US cities are recorded, for example, at \url{http://www.realestate3d.com/gps/uslatlongdegmin.htm}. We convert the longitudes provided by the website above from positive to negative to match with the real data.

\begin{figure}[hbt]
\centering
\includegraphics[width=2.5in]{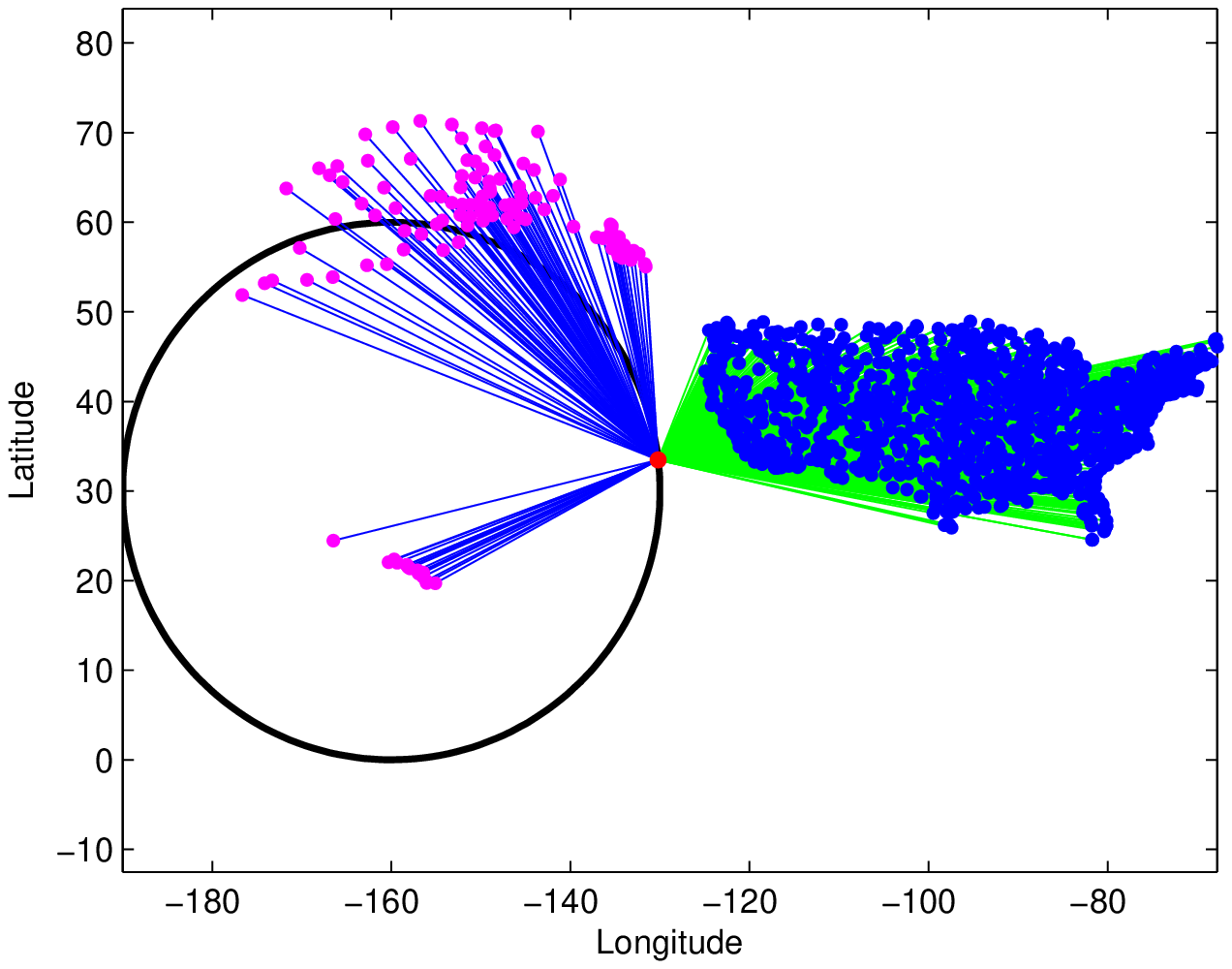}
\includegraphics[width=3.0in]{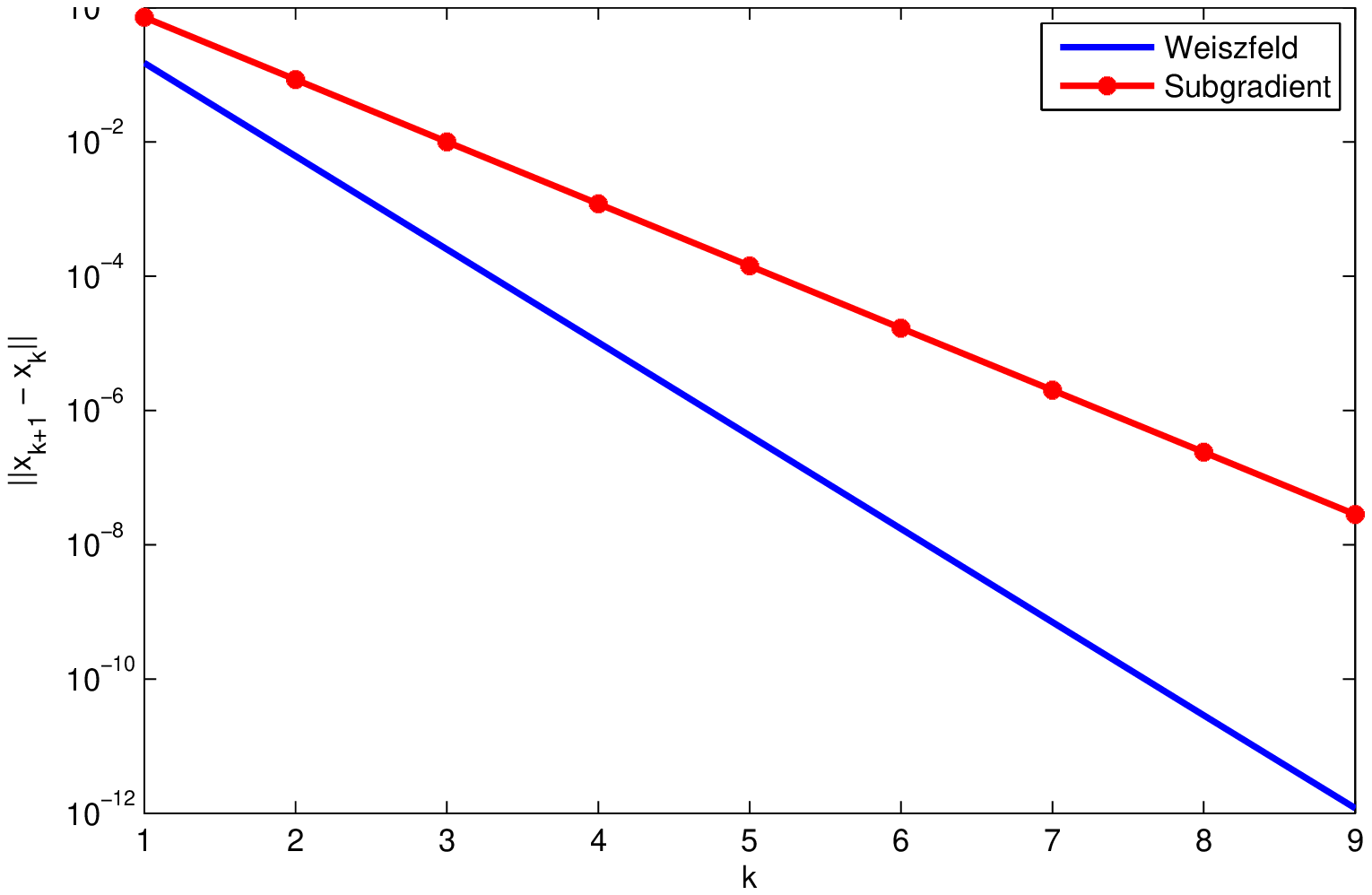}
\caption{A Generalized Fermat-Torricelli Problem with US Cities}
\end{figure}

 We divide these cities into two groups. Group B contains 120 cities in Hawaii and Alaska and the coordinates of each city corresponds to $b_j$ for $j=1,\ldots,120$. Group A contains 1097 other cities and the coordinate of each city corresponds to $a_i$ for $i=1,\ldots, 1097$. We consider problem \eqref{mainproblem} with $\O_i=\{a_i\}$, $\Theta_j=\{b_j\}$, and $\alpha_i=1$, $\beta_j=1$ for all $i,j$. The constraint set is the ball with center at (30, -160) and radius $r=30$.

\medskip
 At each iteration, to find $x_{k+1}$, we solve problem $(P_{y_k})$ by two methods: the generalized Weiszfeld algorithm and the subgradient algorithm (step size $\alpha_\ell=1/\ell$) with the starting point $z_k:=x_k$ and 10 iterations. Both methods give an approximate optimal solution of $x^*\approx(33.48, -130.20)$ and the approximate optimal value $V^*\approx  34,212.80$ although the DCA with generalized Weiszfeld algorithm has much better convergence rate. The subgradient method, however, is applicable to the case where the constraint set $S$ intersects some target sets $\O_i$ for $i\in I$.

\medskip
  The algorithm presented also shows its effectiveness to deal with  set models. Consider problem \eqref{mainproblem} generated by a collection of square target sets belonging to two groups as follows. Group B contains 120 squares whose centers are points in the group B above. Similarly,  Group A contains 1097 squares with centers being points in the group A above.  All squares in both group A and group B have radius of 5 (half-side length). The constraint set is the ball with center $(30,-160)$ and radius 30. We also set $\alpha_i=\beta_j=1$ for all $i,j$. Our algorithm yields an approximate optimal solution of $x^*\approx (32.29, -130.09)$ and an approximate of $V^*\approx  28,803.24$

\medskip
 After getting the approximate optimal solution, we create several randomly selected points in the constraint set and compare the function value with the approximate function value obtained from the DCA. The numerical results experimentally show that the algorithm yields an approximate global optimal solution, although there is no proof that the algorithm converges to a global optimal solution.
 }\end{Example}

\end{document}